\documentclass{article}
\usepackage{graphicx}
\usepackage{amsmath}
\usepackage{amssymb}

\usepackage{amsthm}
\usepackage{makeidx}

\usepackage{color}  
\usepackage{tikz}

\usepackage{hyperref}

\newcommand{\R}{\mathbb{R}}

\newcommand{\cx}{\check{X}}

\usepackage{latexsym}
\usepackage{epsfig}

\usepackage{array}
\usepackage{enumerate}
\newcommand{\bel}[1]{\begin{equation}\label{#1}}

\newcommand{\be}{\begin{equation}}
\newcommand{\ee}{\end{equation}}
\newcommand{\ba}{\begin{eqnarray}}
\newcommand{\ea}{\end{eqnarray}}
\newcommand{\rf}[1]{(\ref{#1})}

\newcommand{\qe}{\end{equation}}

\theoremstyle{theorem}

\theoremstyle{corollary}
\newtheorem{coro}{Corollary}[section]
\theoremstyle{lemma}
\newtheorem{lemma}{Lemma}[section]
\theoremstyle{definition}
\newtheorem{defi}{Definition}[section]
\theoremstyle{remark}
\newtheorem*{pf}{Proof}
\theoremstyle{remark}
\newtheorem*{rem}{Remark}
\theoremstyle{remark}

\title{Topological representation of the geometry of metric spaces}
\author{Parvaneh Joharinad\thanks{p.joharinad@iasbs.ac.ir, Institute for Advanced Studies in Basic Sciences, Zanjan, Iran},   J\"urgen Jost\thanks{jost@mis.mpg.de,  
Max Planck Institute for Mathematics in the Sciences, Leipzig, Germany}}
\begin{document}

\maketitle
\begin{abstract}
 Inspired by the concept of hyperconvexity and its relation to curvature, we  translate geometric properties of a metric space encoded by curvature inequalities into the persistent homology induced by the \v{C}ech filtration of that space.
\end{abstract}
\section{Introduction}
The fundamental concept of sectional curvature was introduced by Riemann \cite{Riemann:1868} for spaces that in his honor are called Riemannian manifolds. The model space is Euclidean space which has vanishing curvature, and for a general Riemannian manifold, its sectional curvature measures its local deviation from being Euclidean. Spheres have positive curvature, whereas hyperbolic spaces have negative curvature.  For more general metric spaces, one can only define curvature inequalities. These reduce to those holding for sectional curvatures in the case of Riemannian manifolds, and they encode important local and global properties of the space in question. In particular, we have the concepts of Alexandrov \cite{Alexandrov:1957} and Busemann \cite{Busemann:1955}; see also \cite{Berestovskij:Nikolaev:1993, Jost:1997} for systematic treatments. Gromov \cite{Gromov:1999} developed a general theory of hyperbolic spaces, that is, spaces whose curvature is negative. Abstract versions of curvature inequalities involve relations between the distances in certain configurations of 4 points, see for instance \cite{Alexander:2019,Berg and Nikolaev:2018,Jost:et:al: 2015}. \\
Another important theory, developed by Isbell \cite{Isbell:1964} and Dress \cite{Dress:1984} and others, e.g.  \cite{Dress et al:2002,Lang:2013}, associates to each metric space a hyperconvex space containing it; the notion of hyperconvexity had been introduced in \cite{Aronszajn:Panitchpakdi:1956}. \\
The concepts of metric curvature inequalities as cited above would naturally associate the curvature $-\infty$ to such hyperconvex spaces. That is, they represent extreme cases for curvature. In \cite{JJ1}, we have therefore introduced an approach to metric curvature for which hyperconvex spaces, instead of Euclidean ones, constitute the natural model spaces. We also pointed to a link with the topological concept of persistent homology \cite{Carlsson:2009}. \\
In the present paper, we shall introduce a systematic construction that translates the geometric properties of a metric space $(X,d)$ encoded by curvature inequalities into the homology of a simplicial complex  $\cx$ associated to that space. That is, we systematically translate local (and global) geometry into topology. Among other things, this also makes the link with persistent homology alluded to above more transparent.  \\
The relation between hyperconvexity and persistent homology is also explored in  \cite{Lim:Memoli:Okutan:2020}, where  the persistent homology of the Vietoris-Rips simplicial filtration is studied using the persistent homology corresponding  to the \v{C}ech complex of a specific covering of $r$- neighbourhoods ($r>0$) of  the metric space in its injective envelope. \\
\noindent
{\bf Some notation:} Let $(X,d)$ be a metric space. We denote by 
\bel{1}
B(x,r):=\{y\in X: d(x,y)\le r\}
\qe
the closed ball with center $x\in X$ and radius $r\ge 0$.

\section{Nerve of a cover and  \v{C}ech homology}

We first recall the notion of a simplicial complex.
\begin{defi}
  A {\it simplicial complex} $\Sigma$ consists of a vertex set $I$ whose members can span simplices that satisfy the following condition. Whenever $I'\subset I$ spans a simplex, then also every $I''\subset I'$ spans a simplex. \\
Equivalently, a simplicial complex $\Sigma$ with vertex set $I$ is a collection of subsets of $I$ with the condition that whenever $I'\in \Sigma$, then also $I''\in \Sigma$ for every  $I''\subset I'$.\\
A simplex with $(p+1)$ vertices is called a $p$-simplex. 
When $S$ is a simplex in the simplicial complex $\Sigma$ with vertex set $I'$, then for $i\in I'$, the simplex spanned by $I\setminus \{i\}$ is called a {\it facet} of $\Sigma$. 
\end{defi}
Thus, a facet of a $p$-simplex is a $(p-1)$-simplex.
We next recall how to construct a simplicial complex from a cover of a topological space.
\begin{defi}
  Let $\mathfrak{U}=\{\mathcal{U}_{\alpha}\}_{\alpha\in I}$ be a cover of the topological space $X$. The {\it nerve} $ \mathcal{N}(\mathfrak{U})$ of $\mathfrak{U}$ is the  simplicial complex with vertex set $I$, where a subset $\{\alpha_{i_0},\alpha_{i_1 },...,\alpha_{i_p}\}\subset I$ defines a $p$ simplex  in $ \mathcal{N}(\mathfrak{U})$ iff 
\be
\bigcap_{k=0}^{k=p}\mathcal{U}_{\alpha_{i_k}}\neq\emptyset
\qe 
\end{defi}

The intention then is to choose a  cover $\mathfrak{U}$ that is good in the sense that its nerve is a good combinatorial model for $X$. 
In fact, by the Nerve Theorem, if $F$ is a (possibly infinite) collection of closed sets in a Euclidean space, where all the non-empty intersections of subcollections are contractible and the union is triangulizable, then the nerve of $F$ has the same homotopy type as the union of its members.\\
If we have a data set  $\mathbb{X}$  in the Euclidean space $\mathbb{E}^d$  and $\mathfrak{U}=\{B(x_{\alpha}, \epsilon), \alpha\in I, x_{\alpha}\in F\subset \mathbb{X}, \epsilon>0\}$ is a cover of closed balls  for $\mathbb{X}$, then the nerve of $\mathfrak{U}$ is topologically equivalent to the $\epsilon$-neighborhood of $F$ in $\mathbb{X}$.\footnote{Here, one can define the distance balls $B(x, \epsilon) =\{ y\in \mathbb{X}: \mathrm{dist}(x,y)\le \epsilon\}$ either with respect to the intrinsic metric of $\mathbb{X}$ or that of the ambient Euclidean space $\mathbb{E}^d$. In the persistent homology theory of Carlsson, either seems possible.  Another approach, used by Ghrist, Edelsbrunner and  Harer, considers  intersections in the ambient space $\mathbb{E}^d$  of closed balls  whose centers lie on $\mathbb{X}$.}\\
The concept of a nerve was introduced by Alexandroff \cite{Alex}, and nerve theorems for sufficiently good (e.g. paracompact) topological spaces were proved by Borsuk \cite{Borsuk}, Leray \cite{Leray} and Weil \cite{Weil}.\\

Importantly, the topology or geometry of the space $X$ imposes restrictions on the combinatorial properties of coverings. 
Let $F$ be a collection of convex closed sets in $\R^d$. Then by Helly's theorem, all the sets in $F$ have a non-empty common intersection  if and only if every subcollection of $F$ with $d+1$ members has a non-void intersection. This theorem imposes a restriction on the structure of the nerve, in the way that if a subcollection $F$ with $k\geq d+1$ members of the cover $\mathfrak{U}=\{\mathcal{U}_{\alpha}\}_{\alpha\in I}$  has all its $d$-simplices in the nerve, then this subcollection itself must be a $k$-simplex in the nerve  $\mathcal{N}(\mathfrak{U})$.\\
This theorem motivates the definition of hyperconvex and  $m$-hyperconvex metric spaces. We recall this concept, first introduced by Aronszajn and Panitchpakdi \cite{Aronszajn:Panitchpakdi:1956}, here.    
\begin{defi}\label{hyperconvex}
The metric space $(X,d)$  is  {\em hyperconvex} if for any family $\{x_i\}_{i\in I}\subset X$ and positive numbers $\{r_i\}_{i\in I}$ with 
\bel{2a}
r_i+r_j\geq d(x_i,x_j) \text{  for }i,j\in I,
\qe
we have
\[
\bigcap_{i\in I}B(x_i,r_i)\neq \emptyset. 
\]
When this holds for all families with at most $m-1$ members, the space is called $m$-hyperconvex 
\end{defi}
\begin{rem}
  Hyperconvexity refines Menger's convexity condition \cite{Menger:1928}, which requires that whenever
\bel{2b}
r_1+r_2\geq d(x_1,x_2),
\qe
we have 
\bel{2c}
B(x_1,r_1)\cap B(x_2,r_2) \neq \emptyset.
\qe
\end{rem}
If $X$ is a hyperconvex (resp. $m$-hyperconvex) space and $\mathfrak{U}=\{B(x_i,r_i)\}_{i\in I}$ is a cover of closed balls, the nerve of this cover satisfies the following property. 
For any subset $J$ (resp. any subset $J$ with at most $m-1$ members) of $I$, the nerve $\mathcal{N}(\mathfrak{U})$ carries a simplex spanned by $J$ if and only if it carries all 1 dimensional faces, i.e. edges, of this simplex.\\
Hyperconvex spaces are special cases of more general spaces called $\lambda$-hyperconvex spaces.
\begin{defi}
A metric space $(X,d)$ is said to be $\lambda$-hyperconvex for $\lambda\geq1$ if for every family $\{B(x_{i}, r_{i})\}_{i\in I}$ of closed balls in $X$ with the property $r_{i}+ r_{j}\geq d(x_{i}, x_{j})$, one has
\bel{del2}\bigcap\limits_{i\in I}B(x_{i},\lambda r_{i})\neq\emptyset.
\qe
\end{defi}
Clearly, there is an issue with the terminology here, as we need to distinguish between the conditions of  $m$-hyperconvexity and $\lambda$-hyperconvexity. Latin letters are used to indicate the number of intersecting closed balls, while Greek letters are used for the scale of expansion to get an intersection of a family of closed balls. \\

In \cite{JJ1}, we defined a function that assigns to each triple of non-collinear points $(x_1,x_2,x_3)$ in the metric space the quantity
\be \label{rho}
\rho(x_1,x_2,x_3):=\inf_{x\in X}\max_{i=1,2,3}\frac{d(x_i,x)}{r_i},
\qe
where $r_1,r_2,r_3$ are obtained by the Gromov products via 
\begin{align}\label{Gromov product}
\nonumber r_1&=\frac{1}{2}(d(x_1,x_2)+d(x_1,x_3)-d(x_2,x_3)),&\\
\nonumber r_2&=\frac{1}{2}(d(x_1,x_2)+d(x_2,x_3)-d(x_1,x_3)),&\\
r_3&=\frac{1}{2}(d(x_1,x_3)+d(x_2,x_3)-d(x_1,x_2)).&
\end{align} 
In fact, let $\mathfrak{U}=\{B(x_i,r_i)\}_{i\in I}$ be a cover of closed balls. If all the edges of $\{i,j,k\}\subset I$ belong to $ \mathcal{N}(\mathfrak{U})$, then $\{i,j,k\}$ spans a simplex in $ \mathcal{N}(\tilde{\mathfrak{U}})$, where  $\tilde{\mathfrak{U}}=\{B(x_i,\rho r_i)\}_{i\in I}$ and 
\[
\rho:=\sup_{i_1,i_2,i_3\in I}\rho(x_{i_1},x_{i_2},x_{i_3}).
\]
The extreme case where the function $\rho$ in (\ref{rho}) is exact (i.e. the infimum is attained) and equals the minimal possible  value $1$ for each triple of non-collinear points is used as a definition of  tripod spaces in \cite{JJ1} or $4$-hyperconvex spaces in \cite{Aronszajn:Panitchpakdi:1956}. The program of \cite{JJ1} is to systematically compare the topology of a metric space with that of such a tripod space.\\

\subsection*{Homology}
\begin{defi}
The  homology groups of the simplicial complex $ \mathcal{N}(\mathfrak{U})$  are called the \v{C}ech homology groups of $(X,d)$ corresponding to the cover $\mathfrak{U}$ and  denoted by $\check{H}_p(\mathfrak{U})$.
\end{defi}
In order to obtain invariants of $(X,d)$ that do not depend on the choice of a cover, one then takes a directed limit of the homology groups of the nerves over all possible covers of $X$, ordered by refinement. While such a construction is possible in principle for arbitrary topological spaces, here we shall present  the definition for metric spaces, as these are the spaces that we are concerned with in topological and geometric data analysis.
\begin{defi}
A pre-ordered set $\mathcal{P}$ is said to be directed, if for every $r,s\in\mathcal{P}$ there exists some $t\in\mathcal{P}$ such that $r\leq t$ and $s\leq t$.
Let $\mathcal{P}'$ be a directed set. An inverse system of objects in the category $\underline{C}$ directed by $\mathcal{P}'$ is a family $\{c_r\}_{r\in \mathcal{P}'}$ of objects of $\underline{C}$ together with morphisms (called projections) $\pi_{sr}:c_s\longrightarrow c_r$ whenever $r\leq s$, such that
\begin{enumerate}
\item[a)]  $\pi_{rr}$ is the identity map on $c_r$
\item[b)] $\pi_{tr}=\pi_{ts}\circ\pi_{sr}$,
whenever $r\leq s\leq t$.
\end{enumerate}
\end{defi}
Given an inverse system $\{c_r, \pi_{rs}\}$ directed by $\mathcal{P}$, by the inverse limit of this system one means the subset $c_{\infty}$ of the Cartesian product $C:=\Pi_{r\in\mathcal{P}}c_r$ defined by
\[
c_{\infty}:=\{(\sigma_r)_{r\in\mathcal{P}}\in C:\: \pi_{ts}(\sigma_t)=\sigma_s\}
\]
Set $\mathcal{P}$ to be the directed set of nonnegative functions $r:X\longrightarrow \R $ and $\mathcal{P}'$ the directed set of  covers $\mathfrak{U}_r:=\{B(x,r(x))\}_{x\in X}$ with $r\in\mathcal{P}$ and the pre-order defined by refinement. Let $ \mathcal{N}(\mathfrak{U}_r)$ be the nerve of the cover $\mathfrak{U}_r$ for each $r\in\mathcal{P}$. One can see that $\mathfrak{U}_r$ is a refinement of $\mathfrak{U}_{r'}$ (i.e. $\mathfrak{U}_{r'}\leq \mathfrak{U}_r$)  whenever $r\leq r'$ (i.e. $r(x)\leq r'(x)$ for all $x\in X$) and consequently $ \mathcal{N}(\mathfrak{U}_r)\subseteq\mathcal{N}(\mathfrak{U}_{r'})$.\\
Therefore, $\{\mathcal{N}(\mathfrak{U}_r)\}_{r\in\mathcal{P}}$ (resp. $\{H_{*}(\mathcal{N}(\mathfrak{U}_r))\}_{r\in\mathcal{P}}$)  defines an inverse system of objects in the category $\underline{C}$ of all simplicial complexes with vertex set $X$  (resp. $F$-vector spaces) directed by $\mathcal{P}'$, where the projection $\pi_{r'r}$ is the inclusion map $\iota:\mathcal{N}(\mathfrak{U}_r)\longrightarrow\mathcal{N}(\mathfrak{U}_{r'})$ (resp. the homomorphism induced by inclusion) for $r\leq r'$.  
\begin{defi}
Let $X$, be a topological space. The \v{C}ech homology group $\check{H}_k(X)$ is defined as the directed limit of $\{H_{k}(\mathcal{N}(\mathfrak{U}))\}_{\mathfrak{U}\in Cov(X)}$, where $Cov(X)$ is the directed set of all coverings of $X$ with the order defined by refinement. 
\end{defi}
In the special case where $(X,d)$ is a metric space, $(\sigma_r)_{r\in\mathcal{P}}$ belongs to $\check{H}_k(X)$, if for each $r\in\mathcal{P}$, one has $\sigma_r\in H_k(\mathcal{N}(\mathfrak{U}_r))$ with $\sigma_{r'}=\iota_*(\sigma_r)$ whenever $r\leq r'$.
By elementary simplicial collapses, one can modify $\mathcal{N}(\mathfrak{U}_r)$ such that each $(k-1)$-cycle consists of $k+1$ simplices (i.e is the boundary of a single empty $k$-simplex) without changing the homology groups. Therefore,  \v{C}ech homology detects nontrivial topology when from a collection of $k$ balls, any $k-1$ of them have a nontrivial  intersection, while the intersection of all $k$ is empty.\\
 We observe the following result which is a direct consequence of the definition of hyperconvexity.
\begin{lemma}
Let $X$ be a hyperconvex space. For each function $r\in\mathcal{P}$ with $r(x)+r(y)\geq d(x,y)$, the associated homology $H_*(\mathfrak{U}_{r'})$ is trivial for each $r'\geq r$.
\end{lemma}
\begin{defi}
 A function $r$ that is  minimal subject to $r(x)+r(y)\geq d(x,y)$ is called an extremal function.
\end{defi}
We omit the proof of the following easy lemma. 
 \begin{lemma}
\begin{enumerate}
\item Every extremal function is $1-$Lipschitz. 
\item For each $x\in X$, the distance function $d(x,.)$ is extremal.
\item If $X$ is compact and $r$ is extremal, then for each $x\in X$ there exists some $y\in X$ with $r(x)+r(y)=d(x,y)$.
\item For each function $r':X\longrightarrow \R $ with $r'(x)+r'(y)\geq d(x,y)$, there exists an extremal function $r\leq r'$.
\end{enumerate}
 \end{lemma}
 
\section{Persistent homology}
Let $(X,d)$ be a connected metric space. 
\begin{defi}\label{sc}
We construct a simplicial complex $\check{X}$ whose vertex set is $\{(x,r): x\in X, r\ge 0\}$. For every index set $I$ and any family $(x_i,r_i)_{i\in I}$ of vertices where all the $x_i$ are different from each other, it carries a simplex spanned by these vertices whenever
\bel{2}
\bigcap_{i\in I} B(x_i,r_i) \neq \emptyset.
\qe
\end{defi}
On the vertex set of $\check{X}$, we have the topology of $X\times \R^+$, and we say that a family of simplices converges to a limit when the corresponding vertex sets converge to one that spans a simplex which then is that limit. Since some of the $x_i$ may coincide in the limit, the limiting simplex could be of smaller dimension than the approximating ones. 
When we consider each simplex as a geometric simplex, that is, the unit simplex in $\R^I$ with its Euclidean topology, then $\cx$ becomes a topological space. 
When $I$ is finite, with $k+1$ elements, we also call the resulting simplex a $k$-simplex.
\begin{rem}
  \begin{itemize}
  \item 
As mentioned in \cite{Carlsson:2009}, this construction is computationally expensive since it requires the storage of simplices of various dimensions.
\item As should become clear in the sequel, the construction does not yield a natural metric on the simplicial complex   $\check{X}$. We could, of course, equip each simplex with a Euclidean metric, but then the scale is not determined. 
\end{itemize}
\end{rem}
Instead of the \v{C}ech complex, we can also take the  Vietoris-Rips complex  $VR(X)$,  which has the same vertex set as $\check{X}$ and carries a simplex spanned by vertices $\{(x_i,r_i): i\in I\}$ whenever 
\bel{2}
r_i+r_j\geq d(x_i,x_j) , \; \forall i,j\in I.
\qe
This complex is obtained  solely from the edge information.  Clearly, both constructions have the same set of edges and $\cx\subset VR(X)$. In fact, the Vietoris-Rips complex is a flag complex, i.e.  it is maximal among all simplicial complexes with the same 1-skeleton. If $X$ is hyperconvex, then the above complexes are identical.
\begin{lemma}\label{lem1}
  As a topological space, $\cx$ is connected.
\end{lemma}
\begin{pf}
  For each $x\in X$, the halfline $\{(x,r): r\ge 0\}$ is connected. And for $x_1,x_2\in X$, $(x_1,r_1)$ and $(x_2,r_2)$ are connected by a 1-simplex as soon as $r_1+r_2\ge d(x_1,x_2)$. And since we assume that $X$ is connected, the distance between any two points is finite. \qed
\end{pf}

\begin{defi}
  For a function $r: X\to \R^+$, the space $\{(x,r(x)\}$ equipped with the simplices according to Def. \ref{sc} is called a {\em slice} of $\cx$. We shall denote it by $\cx_{r(x)}$ or also simply by $\cx_r$. \\
Slices $VR(X,r)$ of $VR(X)$ are defined in the same manner.
  \end{defi}
  \begin{lemma}
 The \v{C}ech homology  $\check{H}_*(X)$ is the directed limit of $\{H_*(\cx_r)\}_{r\in\mathcal{P}}$, where $\mathcal{P}$ is the directed set of all non-negative functions on $X$.
 \end{lemma}

\qed
\begin{lemma}
 For each function $r$, we have $\cx_r\subseteq VR(X,r)$. Moreover, $VR(X,r)\subseteq \cx(\mu r)$, where $\mu>$  expansion constant of $X$. If the expansion constant is exact, one can choose $\mu =$ expansion constant of $X$.
\end{lemma}
\qed

A  function $r$ chooses a point on each fiber $\{x\}\times \R^+$. Then we construct the simplicial complex $\cx_r$ according to  Def. \ref{sc}. If $\cx_r$ contains an $m$-simplex whose vertices lie on fibers over $x_0,x_1,..., x_m$, we have 
\[
r(x_i)+r(x_j)\geq d(x_i,x_j), \; \text{for all}\; 0\leq i,j\leq m.
\]
The converse does not necessarily hold. However, if we construct Vietoris-Rips simplices with respect to the function $r$, the converse is also true. \\
For each subset of $m+1$ points $\{x_0,x_1,...,x_m\}$ in X, let $r:X\longrightarrow \R $ be a function whose restriction to $\{x_0,x_1,...,x_m\}$ is minimal subjected to $r(x_i)+r(x_j)\geq d(x_i,x_j)$ for all $0\leq i<j\leq m$, i.e. is an extremal function on $\{x_0,x_1,...,x_m\}$. (In fact,  such a function is an extension of an element of the hyperconvex envelope of $\{x_0,x_1,...,x_m\}$, which by the way is a simplicial complex, to the whole space $X$.)\\
Then $r$ can be considered as the emerging time of this simplex in the  Vietoris-Rips complex. More precisely, the slice $\cx_r'$ does not have an $m$-simplex with vertex set $\{x_0,x_1,...,x_m\}$ (with an abuse of notation) if $r'\leq r$ on the set $\{x_0,x_1,...,x_m\}$ with strict inequality for at least one point since at least one of the edges is missing in this case. 
\begin{lemma}
If $\cx_{r'}$ contains all the boundary facets of a simplex with vertex set $\{(x_i,r'(x_i))\}_{i\in I}$ with $|I|\ge 3$, then there is a function $r\leq r'$ which is extremal on $\{x_i\}_{i\in I}$. Moreover, if  $X$ satisfies Menger's convexity condition, then $\cx_r$ contains all the edges between $x_i$s for the first time. i.e. if $r''\leq r$ on $\{x_i\}_{i\in I}$ and $\cx_{r''}$ contains all the edges between $x_i$s, then $r''= r$ on $\{x_i\}_{i\in I}$. 
\end{lemma}
Let us now explain how persistent homology \cite{Carlsson:2009} can be seen in this framework. 
  Persistent homology \cite{Carlsson:2009} considers the slices $\cx_r$ for the constant functions $x\mapsto r$ and their homology as $r>0$ varies. 
Now, let $\mathcal{P}$ be the partially ordered set of nonnegative functions $r:X\longrightarrow \R$.  We observe that both constructions (\v{C}ech and Vietoris-Rips) 
were obtained as an inverse system of objects directed by $\mathcal{P}'$ and yield a $\mathcal{P}$-persistent simplicial  complex attached to $X$. \\
In general a $\mathcal{P}$-persistent object in the category $\underline{C}$, where $\mathcal{P}$ is a partially ordered set, is a family $\{c_r\}_{r\in \mathcal{P}}$ of objects of $\underline{C}$ together with morphisms $\phi_{rs}:c_r\longrightarrow c_s$ whenever $r\leq s$, such that
\[
\phi_{rs}=\phi_{rt}\circ\phi_{ts},
\]
whenever $r\leq t\leq s$.\\
The category of all $\mathcal{P}$-persistent objects in $\underline{C}$ is denoted by $\mathcal{P}_{pers}(\underline{C})$. The machinery used in the standard persistent homology method constructs $\R$-persistent simplicial complexes with the vertex set $X$ to study its associated object in $\R_{pers}(\underline{Ab})$, i.e.,  $\R$-persistent homology groups.\\
In fact, let $r\in\R^+$ be the smallest value for which the cycle $\sigma_r$ appears in $H_{k}(\cx_r)$ and $(\sigma_s)_{s\in\mathcal{P}}\in\check{H}_k(X)$ is an element whose projection on $H_{k}(\cx_r)$ is $\sigma_r$. Then the interval $[r,r']$ defines the line corresponding to $\sigma_r$ in the bar code, where $r'>r$ is the smallest value with $proj_{r'}(\sigma_s)_{s\in\mathcal{P}}=0$. \\


\begin{lemma}
 If $X$ is a complete convex metric space, and $r:X\longrightarrow \R$ is an extremal function with respect to the subset $\{x_0,x_1,...,x_m\}$, then there is a number $1\leq \rho\leq 2$ such that $\{x_0,x_1,...,x_m\}$ forms a simplex in the slice $\cx_{\rho r}$ 
\end{lemma}
\begin{pf}
By Menger's convexity we have $B(x_i,r(x_i))\cap B(x_j,r(x_j))\neq \emptyset$ for all $0\leq i,j\leq m$. Then for some number $\rho$, which is $\leq$ expansion constant of $X$, one has
\[
\bigcap_{0\leq i\leq m}B(x_i,\rho r(x_i))\neq \emptyset.
\]
\end{pf}

\begin{lemma}\label{lem2}
 If $(X,d)$  is  hyperconvex, then whenever $\cx$ contains all the boundary facets of a simplex with vertex set $\{(x_i,r_i)\}_{i\in I}$ with $|I|\ge 3$, then it also contains the simplex with that vertex set. In particular, for a hyperconvex space, the \v{C}ech homology of $\cx$ is trivial.
\end{lemma}
\begin{pf}
  If $B(y_1,r_1)\cap B(y_2,r_2)\neq \emptyset$, then necessarily $r_1+r_2\ge d(y_1,y_2)$. Thus, when $B(x_i,r_i)\cap B(x_j,r_j)\neq \emptyset$ for all $i,j\in I$, that is, when $\cx$ contains all the edges (1-simplices) between $(x_i,r_i)$ and $(x_j,r_j)$, then $r_i+r_j\geq d(x_i,x_j)$, and therefore, by hyperconvexity, it also contains all the simplices for all vertex subsets of $\{(x_i,r_i)\}_{i\in I}$. This implies the first claim. Therefore, all \v{C}ech cohomology in dimension $\ge1$ is trivial. Since $\cx$ is connected by Lemma \ref{lem1}, this concludes the proof. \qed
\end{pf}
\begin{coro}\label{cor1}
  When $(X,d)$  is  hyperconvex, then every slice $\cx_r$  has trivial homology except possibly in dimension 0 (it may be disconnected) 
\end{coro}
\qed 

\begin{coro}
If $(X,d)$  is $m$-hyperconvex, then the \v{C}ech homology of $\cx$ (resp. every slice $\cx_r$) in dimension $\leq m-2$( resp. in $1\leq \text{dimension}\leq m-2$) is trivial. 
\end{coro}
\qed

In order to put Lemma \ref{lem2} into perspective, we now turn to another example, the Euclidean plane, that is, $\R^2$ equipped with its Euclidean metric. Let $x_1,x_2,x_3$ be equidistant, that is, $d(x_1,x_2)=d(x_2,x_3)=d(x_3,x_1)=:a$. Then
\begin{eqnarray}
\nonumber
B(x_i,a\rho)\cap B(x_j,a\rho)&\neq \emptyset &\text{ iff } \rho \ge \frac{1}{2}\\
  \label{eucl}
  B(x_1,a\rho)\cap B(x_2,a\rho)\cap B(x_3,a\rho) &\neq \emptyset &\text{ iff } \rho \ge \frac{1}{\sqrt{3}},
\end{eqnarray}
and consequently, for $\frac{1}{2} \le \rho < \frac{1}{\sqrt{3}}$, $\check{\R}^2$ contains three edges between the $(x_i,a\rho)$, but not the 2-simplex with those edges.\\
More generally, for three arbitrary points $x_1,x_2,x_3\in \R^2$,  $\check{\R}^2$ contains three edges between the vertices $\{(x_i,r_i),\; i=1,2,3\}$, where $r_1,r_2, r_3$ are the Gromov products obtained by (\ref{Gromov product}),  
\begin{eqnarray}
\nonumber
B(x_i,\rho r_i)\cap B(x_j,\rho r_j)&\neq \emptyset &\text{ iff } \rho \ge 1\\
  \label{eucl}
  B(x_1,\rho r_1)\cap B(x_2,\rho r_2)\cap B(x_3,\rho r_3) &\neq \emptyset &\text{ iff } \rho \ge \rho(x_1,x_2,x_3),
  \end{eqnarray}
  where $\rho(x_1,x_2,x_3)$ is obtained by (\ref{rho}) and is a number strictly larger than $1$ with maximal value $\frac{2}{\sqrt{3}}$. For  $\rho=\rho(x_1,x_2,x_3)$ the intersection of three scaled balls is a single point, called weighted circumcenter, whose distance from each $x_i$ is equal to $\rho r_i$. \\
 For a general metric space $(X,d)$, \cite{Jost:et:al: 2015} would suggest to compare the length of the interval of those $r
\rho >0$ for which we find unfilled simplices spanned by three equidistant vertices  with the length $\frac{1}{\sqrt{3}} -\frac{1}{2}$ of that interval for the Euclidean plane and to take that as a measure of curvature. The curvature is negative (positive) if that length is smaller (larger) than the Euclidean one. In \cite{JJ1}, however, we have suggested to take those spaces where the length of that interval vanishes, the {\em tripod spaces}, as the comparison spaces, instead of the Euclidean plane. Hyperconvex spaces are tripod spaces, as follows for instance from Cor. \ref{cor1}.  \\

While in \cite{Jost:et:al: 2015,JJ1}, only configurations of three points have been considered, that is, putting it into the framework developed here, we have only looked at the first homology group of $\cx$, it is now natural to look at the entire homology of $\cx$. That is, we systematically translate the local geometry of a metric space $(X,d)$ into the topology, or more precisely, the homology of the space $\cx$.\\

In fact, $\cx$ not only incorporates the local geometry of $(X,d)$, but also its global topology. To see, let us consider another simple example, the unit circle obtained by identifying the two endpoints $0$ and $1$ of the unit interval $[0,1]$. We consider three equidistant points, for instance $x_1=0, x_2=\frac{1}{3}, x_3=\frac{2}{3}$. Then 
\begin{eqnarray}
\nonumber
B(x_i,r)\cap B(x_j,r)&\neq \emptyset&,\\
\nonumber
 \text{ but } B(x_1,r)\cap B(x_2,r) \cap B(x_3,r) &=\emptyset &\text{ for } \frac{1}{6}\le r < \frac{1}{3}\\
  \label{4}
  B(x_1,\frac{1}{3})\cap B(x_2,\frac{1}{3}) \cap B(x_3,\frac{1}{3})=\{x_1,x_2,x_3\} &\neq \emptyset&.
\end{eqnarray}
Therefore, in the range $\frac{1}{6}\le r < \frac{1}{3}$, our simplicial complex contains three 1-simplices with vertex set $(x_1,r), (x_2,r),(x_3,r)$ that are not filled by a 2-simplex. \\
Note: The three radii obtained by the Gromov product for these three points are $r_1=r_2=r_3=1/6$ and the function $\rho$ reaches its maximal value 2 at $(x_1,x_2,x_3)$. This is the largest possible value, since the expansion constant of a complete metric space is not larger than $2$.

\section{Computations}
For each triple $(x_1,x_2,x_3)$ let $\lambda$ be the supremum of the numbers $\alpha$ that satisfy
\begin{align}
\nonumber\alpha d(x_1,x_2)\leq d(x_1,x_3)+d(x_2,x_3),&\\
\nonumber \alpha d(x_1,x_3)\leq d(x_1,x_2)+d(x_2,x_3),&\\
 \alpha d(x_2,x_3)\leq d(x_1,x_2)+d(x_1,x_3).&
\end{align} 

If the supremum is attained, i.e. $\lambda$ satisfies the above inequalities, then at least one of them has to be an equality. If $\lambda=1$, the three points are collinear, the corresponding triangle is degenerate. If $\lambda=2$, all inequalities become equalities and the correspongind triangle is equilateral. \\
We can classify all triangles in $(X,d)$ with respect to this measure. Let us put
\[
X_{\lambda}:=\{ \text{all triangles with measure equal to $\lambda$}\}.
\]
For instance, $X_1$ is the set of all degenerate triangles and $X_2$ is the set of all equilateral triangles and every other triangle lies somewhere between theses two classes.
\begin{lemma}
Assume for $(x_1,x_2,x_3)\in X_{\lambda}$, equality $\lambda d(x_i,x_j)= d(x_i,x_k)+d(x_j,x_k)$ occurs for a permutation $(i,j,k)$ of $(1,2,3)$. Then the edge $[x_i,x_j]$ has the largest length and $r_k=\dfrac{\lambda-1}{2}d(x_i,x_j)$.
\end{lemma}
\begin{pf}
Assume
\begin{align*}
\lambda d(x_i,x_j)= d(x_i,x_k)+d(x_j,x_k),&\\
 \lambda d(x_i,x_k)\leq d(x_i,x_j)+d(x_j,x_k),&\\
\lambda d(x_j,x_k)\leq d(x_i,x_j)+d(x_i,x_k).&
\end{align*}
Then by direct computation $r_k=\dfrac{\lambda-1}{2}d(x_i,x_j)$. On the other hand, the equality along with the first inequality, resp. second inequality, implies $r_k\leq \frac{1}{2}d(x_i,x_k)$, resp. $r_k\leq \frac{1}{2}d(x_j,x_k)$. Therefore, both $r_j$ and $r_i$ are bigger than $r_k$ since $r_i+r_k=d(x_i,x_k)$ and $r_j+r_k=d(x_j,x_k)$. This  completes the proof. \qed
\end{pf}
According to this lemma, The triangle $(x_1,x_2,x_3)$ is in $X_{\lambda}$, if and only if $\lambda$ is equal to the sum of the lengths of the two smallest sides divided by the length of the largest side.\\    

We can draw the corresponding diagrams for some values of $\lambda$ (e.g. $\lambda =\dfrac{5}{4}, \dfrac{6}{4}, \dfrac{7}{4}, 2$). The diagram corresponding to $\lambda=2$ gives us the comparison of equilateral triangles with their counterparts in the Euclidean plane, the circle and the tree.\\
If $r$ is half  the perimeter of the triangle, then
\bel{200}
\rho(x_1,x_2,x_3)=\frac{3}{r}\min_{x}\max_{i=1,2,3} d(x,x_i)
\qe 
for each triangle $(x_1,x_2,x_3)\in X_2$.\\
Let $(\bar{x}_1,\bar{x}_2,\bar{x}_3)$, $(\bar{\bar x}_1,\bar{\bar x}_2,\bar{\bar x}_3)$ and $(x'_1,x'_2,x'_3)$ be the comparison triangles in the Euclidean plane, the tree and the circle with the same perimeter as the original triangle. Then according to the previous section, 
\begin{align*}
&\rho(\bar{x}_1,\bar{x}_2,\bar{x}_3)=\dfrac{2}{\sqrt{3}},& \\
&\rho(\bar{\bar x}_1,\bar{\bar x}_2,\bar{\bar x}_3)=1,& \\
&\rho(x'_1,x'_2,x'_3)=2, &
\end{align*}
for each triangle $(x_1,x_2,x_3)\in X_2$.\\
For instance, the $(r,\rho)$ diagram for hyperbolic space looks the same as for the Euclidean plane for sufficiently small values of $r$, since smooth manifolds are locally Euclidean, and converges to the line $\rho=1$ at infinity.

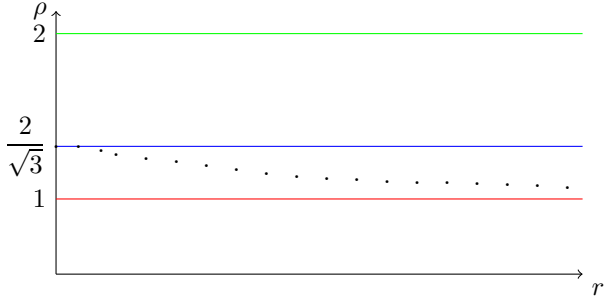
\begin{tikzpicture}[domain=0:2] 
\draw[->] (0,0) -- (7,0) node[below right] {$r$}; 
\draw[->] (0,0) -- (0,3.5) node[left] {$\rho$};
\draw (-0.01,1) node[anchor=east] {1};
\draw  [color=red] (0,1) -- (7,1) ;
\draw (-0.01,1.7) node[anchor=east] {$\dfrac{2}{\sqrt{3}}$};
\draw [color=blue] (0,1.7) -- (7,1.7) ;
\draw (-0.01,3.2) node[anchor=east] {2};
\draw [color=green] (0,3.2) -- (7,3.2) ;
\draw (0,1.56) node[anchor=south] {.};
\draw (0.3, 1.56) node[anchor=south] {.};
\draw (0.6, 1.5) node[anchor=south] {.};
\draw (0.8, 1.45) node[anchor=south] {.};
\draw (1.2, 1.40) node[anchor=south] {.};
\draw (1.6, 1.35) node[anchor=south] {.};
\draw (2, 1.3) node[anchor=south] {.};
\draw (2.4, 1.25) node[anchor=south] {.};
\draw (2.8, 1.2) node[anchor=south] {.};
\draw (3.2, 1.15) node[anchor=south] {.};
\draw (3.6, 1.13) node[anchor=south] {.};
\draw (4, 1.11) node[anchor=south] {.};
\draw (4.4, 1.09) node[anchor=south] {.};
\draw (4.8, 1.08) node[anchor=south] {.};
\draw (5.2, 1.07) node[anchor=south] {.};
\draw (5.6, 1.06) node[anchor=south] {.};
\draw (6,1.05) node[anchor=south] {.};
\draw (6.4, 1.04) node[anchor=south] {.};
\draw (6.8, 1.01) node[anchor=south] {.};
 \end{tikzpicture}

The fact that for small values of $r$, $\rho$ is similar to the Euclidean value holds for all Riemannian manifolds. When the Riemannian sectional curvature is positive (negative), the values of $\rho$ will increase (decrease) as a function of $r$. That $\rho=2$ for the circle reflects a topological rather than a geometric property. (In fact, the necessary condition for a triangle on a circle of radius $s$ to be non-degenerate is that its perimeter be equal to $2\pi s$. Thus for $\lambda\neq 1$, the sets $X_{\lambda}$ only consists of triangles with maximal $r$, i.e. $\pi s$.)
 When $(X,d)$ is a sphere of radius $s$, then the largest equilateral triangle consists of three equidistant points on a great circle. The maximal $r$ is thus $\pi s$, and for $r=\pi s$, we have $\rho = \frac{3}{2}$ in \rf{200}. The sectional curvature of that sphere is $\frac{1}{s^2}$. Thus, the largest value of $\rho$ does not depend on the sectional  curvature, but the value of $r$ where it is reached is a decreasing function of that curvature. Obviously, $\frac{3}{2}$ is smaller than $2$, but we can also reach $\rho =2$ on a simply connected Riemannian manifold. For instance, we can cap off a long cylinder at its ends. Or we can take a spherical surface with a waist, where that waist then is represented by closed geodesic that is homotopic to a point, but locally minimizing. \\
On a Riemannian manifold with varying curvature, there may be different values of $\rho$ for any $r$, and so, the diagram gets more complicated.\\
Of course, we can compute and draw such a profile also for other values of $\lambda$. \\
For instance, for three arbitrary points $x_1,x_2, x_3\in \R^2$, the value of $\rho$ and the position of the weighted circumcenter $p$ are determined by the following system of equations
\begin{align*}
d^{2}(x_1,x_2)&=\rho^{2}(r_{1}^{2}+r_{2}^{2}-2 r_{1} r_{2}\cos(\alpha)),&\\
d^{2}(x_2,x_3)&=\rho^{2}(r_{2}^{2}+r_{3}^{2}-2 r_{2} r_{3}\cos(\beta)),&\\
d^{2}(x_1,x_3)&=\rho^{2}(r_{1}^{2}+r_{3}^{2}+2 r_{1} r_{3}\cos(\alpha+\beta)),&
\end{align*}
where $\alpha$ and $\beta$ are the angles at $p$ opposite to  $[x_1, x_2]$ and $[x_2, x_3]$ respectively.\\
Each triangle $(x_1,x_2, x_3)$ on the circle of radius $s$ centered at $O$ can be represented by its corresponding angles $\angle_{O}(x_i,x_j)$ and consequently the Gromov products constitute the unique solution of the following system of equations
\[
\theta_i+\theta_j=\angle_{O}(x_i,x_j), \; 1\leq i< j\leq 3.
\]
Let us suppose $\theta_1\geq\theta_2\geq\theta_3$. Then $\rho(x_1,x_2, x_3)=\dfrac{2\pi}{\theta_1+\theta_2}-1=\dfrac{2\pi}{\angle_{O}(x_1,x_2)}-1$.

\end{document}